\documentclass{amsart}
\usepackage{amssymb}
\usepackage{graphicx}
\usepackage{amscd}
\usepackage{amsthm,amsmath}

\newtheorem{theorem}{Theorem}[section]

\theoremstyle{proposition}
\newtheorem{proposition}[theorem]{Proposition}

\theoremstyle{remark}
\newtheorem{remark}[theorem]{Remark}

\begin{document}
\title[Norms of positive definite T\"{o}plitz matrices and signal processing]{Asymptotic estimates of the norms of positive definite T\"{o}plitz matrices
and detection of quasi-periodic components of stationary random signals}
\author{Vadim M. Adamyan}
\address[Vadim M. Adamyan]{\textit{Department of Theoretical Physics}, \textit{I.I. Mechnikov Odessa
National University,} \\
\textit{65026 Odessa, Ukraine}}
\email[Vadim M. Adamyan]{vadamyan@paco.net}
\thanks{V. Adamyan was supported by the USA Civil Research and Development
Foundation and the Government of Ukraine (CRDF grant UM1-2567-OD-03).}
\author{Jos\'{e} Luis Iserte}
\address[Jos\'{e} Luis Iserte]{\textit{Department of Nuclear Engineering, Polytechnic University of
Valencia, }\\
\textit{46022 Valencia, Spain}}
\email[Jos\'{e} Luis Iserte]{joisvi@cap.upv.es}
\author{Igor M. Tkachenko}
\address[Igor M. Tkachenko, corresponding author]{\textit{Department of Applied Mathematics Polytechnic University of
Valencia, }\\
\textit{46022 Valencia, Spain }}
\email[Igor M. Tkachenko]{imtk@mat.upv.es}
\date{June 26, 2005}
\keywords{Norms of positive definite T\"{o}plitz matrices, positive integral operators
with difference kernels, signal processing}
\thanks{Authors thank Professor G. Verd\'{u} Mart\'{i}n for illuminating
discussions, which stimulated this work. The financial support of the
Polytechnic University of Valencia is gratefully acknowledged.}
\maketitle

\begin{abstract}
Asymptotic forms of the Hilbert-Scmidt and Hilbert norms of positive
definite T\"{o}plitz matrices $Q_{N}=\left( b(j-k)\right) _{j,k=0}^{N-1}$ as
$N\rightarrow \infty $ are determined. Here $b(j)$ are consequent
trigonometric moments of a generating non-negative mesure $d\sigma (\theta )$
on $\left[ -\pi ,\pi \right] $. It is proven that $\sigma (\theta )$ is
continuous if and only if any of those contributions is $o(N)$. Analogous
criteria are given for positive integral operators with difference kernels.

Obtained results are applied to processing of stationary random signals, in particular,
neutron signals emitted by boiling water nuclear reactors.
\end{abstract}

\section{Introduction}

This paper is motivated by the studies \cite{VGMNPLCRS} devoted to the
problem of detection of hidden unstable components in random neutron signals
measured in boiling water nuclear reactors. We assume that a signal from a
monitored system forms, during a sufficiently long time interval, a
real-valued stationary random process $\xi \left( t\right) ,$ $\,t\in \Bbb{Z}
$, with discrete time, such that the means
\begin{equation*}
\left\langle \xi \left( t\right) \right\rangle =0,\;\left\langle \xi
^{2}\left( t\right) \right\rangle =1\text{ }.
\end{equation*}
The correlation function of such a process
\begin{equation*}
b\left( t\right) :=\left\langle \xi \left( t\right) \xi \left( 0\right)
\right\rangle =\left\langle \xi \left( t+t^{\prime }\right) \xi \left(
t^{\prime }\right) \right\rangle ,\;t,t^{\prime }\in \Bbb{Z}\text{ ,}
\end{equation*}
is a real-valued sequence, which admits the representation
\begin{equation}
b\left( t\right) =\int_{-\pi }^{\pi }\exp \left( it\theta \right) d\sigma
\left( \theta \right) \text{ },  \label{a1}
\end{equation}
where $\sigma \left( \theta \right) $ is a non-decreasing bounded function
on $\left[ -\pi ,\pi \right] $ \cite{R}. By virtue of our assumptions
\begin{equation*}
b\left( 0\right) =\sigma \left( \pi \right) -\sigma \left( -\pi \right) =1%
\text{ },
\end{equation*}
and for any $\theta _{1},$ $\theta _{2}$ such that $0\leq \theta _{1}\leq
\theta _{2}\leq \pi $, we have
\begin{equation}
\sigma \left( \theta _{2}\right) -\sigma \left( \theta _{1}\right) =\sigma
\left( -\theta _{1}\right) -\sigma \left( -\theta _{2}\right) \text{ }.
\label{a2}
\end{equation}
In general, the spectral distribution function $\sigma \left( \theta \right)
$ which determines the process correlation function by (\ref{a1}), can be
split into a sum
\begin{equation}
\sigma \left( \theta \right) =\sigma _{c}\left( \theta \right) +\sigma
_{d}\left( \theta \right)  \label{aa2}
\end{equation}
of a continuous non-decreasing function $\sigma _{c}\left( \theta \right) $
and a non-decreasing step function $\sigma _{d}\left( \theta \right) $, and
such a representation is unique up to constant contributions in $\sigma
_{c}\left( \theta \right) $ and $\sigma _{d}\left( \theta \right) $ \cite{R}%
. Notice that both functions, $\sigma _{c}\left( \theta \right) $ and $%
\sigma _{d}\left( \theta \right) $, satisfy the condition (\ref{a2}).
Actually, the problem, formulated in \cite{VGMNPLCRS}, was to find, in a
real-time operation mode, whether the spectral distribution function of the
random signal $\sigma \left( \theta \right) $ contains or does not contain a
non-trivial component $\sigma _{d}\left( \theta \right) $.We will call here
for brevity the random process (signal) $\xi \left( t\right) $ \textit{stable%
} if its spectral distribution function $\sigma \left( \theta \right) $ is
continuous and \textit{unstable} otherwise. In other words, the process is
unstable if and only if
\begin{equation*}
\sigma _{d}\left( \pi \right) -\sigma _{d}\left( -\pi \right) >0\text{ }.
\end{equation*}
The main task of the present work is to formulate criteria of stability of
the process in terms of its correlation function $b\left( t\right) $. As a
main tool we use the sequence of positive definite T\"{o}plitz matrices $%
Q_{N}=\left( b\left( j-k\right) \right) _{j,k=0}^{N-1}$. This paper is
organized in the following way.

In Section 2 we find the principal asymptotic contribution of the
Hilbert-Schmidt norm $\left\| Q_{N}\right\| _{2}$ with $N\rightarrow \infty $%
, and show that the process is stable if and only if this contribution is $%
o\left( N\right) $.

Section 3 contains a similar criterion, but with the Hilbert norm $\left\|
Q_{N}\right\| $ instead of $\left\| Q_{N}\right\| _{2}$. We prove here that
if $N\rightarrow \infty $, then $\left\| Q_{N}\right\| =m\cdot N+o(N)$,
where $m$ is the maximal jump of $\sigma \left( \theta \right) $.

In Section 4 both criteria are generalized for continuous time processes or
for positive integral operators with difference kernels.

In Section 5 the efficiency of the above stability criteria for signal
processing is discussed and illustrated by application to simulated signals
and real neutron signals emitted by the Forsmark 1\&2 boiling water reactor.

\section{Asymptotic form of the Hilbert-Schmidt norm of truncated
correlation matrices and the stability criterion}

Let us denote by $\left\{ Q_{N}\right\} ,\;N=1,2,\ldots $, the sequence of
T\"{o}plitz matrices $\left( b\left( j-k\right) \right) _{j,k=0}^{N-1}$ and
let $\left\| Q_{N}\right\| _{2}$ be the Hilbert-Schmidt norm of $Q_{N}$:
\begin{equation}
\left\| Q_{N}\right\| _{2}=\left( \sum\limits_{j,k=0}^{N-1}b^{2}\left(
j-k\right) \right) ^{\frac{1}{2}}=\left( Nb^{2}\left( 0\right)
+2\sum\limits_{k=1}^{N-1}\left( N-k\right) b^{2}\left( k\right) \right) ^{%
\frac{1}{2}}\text{ }.  \label{a3}
\end{equation}
Our assumptions imply that
\begin{equation*}
b^{2}\left( k\right) =\left| \int_{-\pi }^{\pi }\exp \left( ik\theta \right)
d\sigma \left( \theta \right) \right| ^{2}\leq b^{2}\left( 0\right) =1\text{
}.
\end{equation*}
Therefore $\left\| Q_{N}\right\| _{2}\leq N$.

\begin{theorem}
\label{t1} The process $\xi \left( t\right) $ is stable if and only if
\begin{equation*}
\underset{N\rightarrow \infty }{\lim }\frac{1}{N}\left\| Q_{N}\right\| _{2}=0%
\text{ }.
\end{equation*}
\end{theorem}

\proof%
Introduce
\begin{equation*}
\Phi _{N}\left( \theta \right) =\int\limits_{-\pi }^{\pi }\frac{\sin ^{2}%
\frac{1}{2}N\left( \theta -\theta ^{\prime }\right) }{N^{2}\sin ^{2}\frac{1}{%
2}\left( \theta -\theta ^{\prime }\right) }d\sigma _{c}\left( \theta
^{\prime }\right)
\end{equation*}

and
\begin{eqnarray*}
\Psi _{N}\left( \theta \right) &=&\int\limits_{-\pi }^{\pi }\frac{\sin ^{2}%
\frac{1}{2}N\left( \theta -\theta ^{\prime }\right) }{N^{2}\sin ^{2}\frac{1}{%
2}\left( \theta -\theta ^{\prime }\right) }d\sigma _{d}\left( \theta
^{\prime }\right) = \\
&&\frac{1}{2}\sum\limits_{\alpha }\left( \frac{\sin ^{2}\frac{1}{2}N\left(
\theta -\theta _{\alpha }\right) }{N^{2}\sin ^{2}\frac{1}{2}\left( \theta
-\theta _{\alpha }\right) }+\frac{\sin ^{2}\frac{1}{2}N\left( \theta +\theta
_{\alpha }\right) }{N^{2}\sin ^{2}\frac{1}{2}\left( \theta +\theta _{\alpha
}\right) }\right) m_{\alpha }\text{ },
\end{eqnarray*}
where $\left\{ \pm \theta _{\alpha }\right\} $ is the set of jump points of $%
\sigma \left( \theta \right) $ or, what is the same, of the points of growth
of $\sigma _{d}\left( \theta \right) ,$%
\begin{equation*}
\frac{1}{2}m_{\alpha }=\sigma \left( \theta _{\alpha }+0\right) -\sigma
\left( \theta _{\alpha }-0\right) =\sigma _{d}\left( -\theta _{\alpha
}+0\right) -\sigma _{d}\left( -\theta _{\alpha }-0\right) \text{ }.
\end{equation*}
Due to (\ref{a1}) and (\ref{a3}), we have
\begin{equation*}
\frac{1}{N^{2}}\left\| Q_{N}\right\| _{2}^{2}=\frac{1}{N^{2}}%
\sum\limits_{j,k=0}^{N-1}b^{2}\left( j-k\right) =
\end{equation*}
\begin{equation}
\frac{1}{N^{2}}\int\limits_{-\pi }^{\pi }\int\limits_{-\pi }^{\pi }\left(
\sum\limits_{j,k=0}^{N-1}\exp \left( i\left( j-k\right) \left( \theta
-\theta ^{\prime }\right) \right) \right) d\sigma \left( \theta ^{\prime
}\right) d\sigma \left( \theta \right) =  \label{a4}
\end{equation}
\begin{equation*}
\int\limits_{-\pi }^{\pi }\Phi _{N}\left( \theta \right) d\sigma _{c}\left(
\theta \right) +2\int\limits_{-\pi }^{\pi }\Phi _{N}\left( \theta \right)
d\sigma _{d}\left( \theta \right) +\int\limits_{-\pi }^{\pi }\Psi _{N}\left(
\theta \right) d\sigma _{d}\left( \theta \right) \text{ }.
\end{equation*}
Let
\begin{equation*}
\omega \left( \delta \right) =\underset{\left| \theta -\theta ^{\prime
}\right| \leq \delta }{\max }\left| \sigma _{c}\left( \theta \right) -\sigma
_{c}\left( \theta ^{\prime }\right) \right| \text{ }.
\end{equation*}
The continuity of $\sigma _{c}\left( \theta \right) $ implies $\omega \left(
\delta \right) \underset{\delta \downarrow 0}{\rightarrow }0$.

Since
\begin{equation}
\frac{\sin ^{2}\frac{1}{2}N\left( \theta -\theta ^{\prime }\right) }{%
N^{2}\sin ^{2}\frac{1}{2}\left( \theta -\theta ^{\prime }\right) }\leq 1
\label{a5}
\end{equation}
and
\begin{equation*}
\left| \sin x\right| \geq \frac{2}{\pi }\left| x\right| \mbox{ },\quad
\,\left| x\right| \leq \frac{\pi }{2}\mbox{ }\left( {\rm mod}2\pi \right)
\text{ },
\end{equation*}
then, for any $0<\delta <\frac{2}{\pi }$, we have
\begin{equation*}
\Phi _{N}\left( \theta \right) =\int\limits_{\left| \theta -\theta ^{\prime
}\right| <\delta \mbox{ }\left( {\rm mod}2\pi \right) }\frac{\sin ^{2}\frac{%
1}{2}N\left( \theta -\theta ^{\prime }\right) }{N^{2}\sin ^{2}\frac{1}{2}%
\left( \theta -\theta ^{\prime }\right) }d\sigma _{c}\left( \theta ^{\prime
}\right) +
\end{equation*}
\begin{equation*}
\int\limits_{\left| \theta -\theta ^{\prime }\right| >\delta \text{ }\left(
{\rm mod}2\pi \right) }\frac{\sin ^{2}\frac{1}{2}N\left( \theta -\theta
^{\prime }\right) }{N^{2}\sin ^{2}\frac{1}{2}\left( \theta -\theta ^{\prime
}\right) }d\sigma _{c}\left( \theta ^{\prime }\right) \leq
\end{equation*}
\begin{equation}
2\omega \left( \delta \right) +\frac{\pi ^{2}}{N^{2}}\left[ \sigma
_{c}\left( \pi \right) -\sigma _{c}\left( -\pi \right) \right] \leq 2\omega
\left( \delta \right) +\frac{\pi ^{2}}{N^{2}\delta ^{2}}\text{ }.
\label{ab5}
\end{equation}
Therefore, by an appropriate choice of $\delta $ and $N$, the first two
integrals on the right-hand part of (\ref{a4}) can be made arbitrarily
small. Hence, these integrals tend to zero as $N\rightarrow \infty $.

Rewrite now the third integral on the right-hand side of (\ref{a4}) in the
form
\begin{equation*}
\int\limits_{-\pi }^{\pi }\Psi _{N}\left( \theta \right) d\sigma _{d}\left(
\theta \right) =\frac{1}{2}\sum\limits_{\alpha }m_{\alpha }^{2}+
\end{equation*}
\begin{equation*}
\frac{1}{2}\sum\limits_{\alpha ^{\prime }\neq \alpha }m_{\alpha }m_{\alpha
^{\prime }}\left( \frac{\sin ^{2}\frac{1}{2}N\left( \theta _{\alpha }-\theta
_{\alpha ^{\prime }}\right) }{\sin ^{2}\frac{1}{2}\left( \theta _{\alpha
}-\theta _{\alpha ^{\prime }}\right) }+\frac{\sin ^{2}\frac{1}{2}N\left(
\theta _{\alpha }+\theta _{\alpha ^{\prime }}\right) }{\sin ^{2}\frac{1}{2}%
\left( \theta _{\alpha }+\theta _{\alpha ^{\prime }}\right) }\right) \text{ }%
,
\end{equation*}
and assume that $\sigma _{d}\left( \theta \right) $ has a finite number of
jumps, then we can conclude, taking into account (\ref{a5}), (\ref{ab5}) and
the inequality
\begin{equation*}
\sum\limits_{\alpha }m_{\alpha }=\sigma _{d}\left( \pi \right) -\sigma
_{d}\left( -\pi \right) \leq \sigma \left( \pi \right) -\sigma \left( -\pi
\right) =1
\end{equation*}
that
\begin{equation}
\underset{N\rightarrow \infty }{\lim }\frac{1}{N^{2}}\left\| Q_{N}\right\|
_{2}^{2}=\frac{1}{2}\sum\limits_{\alpha }m_{\alpha }^{2}\text{ }.  \label{a7}
\end{equation}
Thus,
\begin{equation}
\begin{array}{l}
\underset{\alpha }{\max }m_{\alpha }\leq \sqrt{\frac{1}{2}%
\sum\limits_{\alpha }m_{\alpha }^{2}}=\underset{N\rightarrow \infty }{\lim }%
\frac{1}{N}\left\| Q_{N}\right\| _{2}\leq \\
\sqrt{\underset{\alpha }{\max }m_{\alpha }}\sqrt{\frac{1}{2}%
\sum\limits_{\alpha }m_{\alpha }}\leq \sqrt{\underset{\alpha }{\max }%
m_{\alpha }}\text{ }.
\end{array}
\label{aa7}
\end{equation}
The generalization of (\ref{a7}) and, hence, of (\ref{aa7}) for the case of $%
\sigma _{d}\left( \theta \right) $ having infinitely many jumps can be
obtained by continuity using the standard method employed in the next
section.%
\endproof%

The condition of Theorem \ref{t1} in special cases can be specified.

\begin{proposition}
\label{pr} If the process is stable and its spectral distribution
function $\sigma \left( \theta \right) $ satisfies the H\"{o}lder condition:
\begin{equation}
\left| \sigma \left( \theta _{1}\right) -\sigma \left( \theta _{2}\right)
\right| \leq A\left| \theta _{1}-\theta _{2}\right| ^{\nu }\text{ ,}\;0<\nu
\leq 1\text{,}\;-\pi \leq \theta _{1},\theta _{2}\leq \pi ,  \label{ab7}
\end{equation}
then
\begin{equation*}
\left\| Q_{N}\right\| _{2}\underset{N\rightarrow \infty }{=}O\left( N^{\frac{%
2}{2+\nu }}\right) \text{ }.
\end{equation*}
\end{proposition}

\proof%
Using the estimate (\ref{ab5}) and the assumptions of the proposition we
deduce that
\begin{equation*}
\left\| Q_{N}\right\| _{2}\leq \sqrt{2A\delta ^{\nu }N^{2}+\frac{\pi ^{2}}{%
\delta ^{2}}},\;0<\delta \leq \frac{\pi }{2}\text{ }.
\end{equation*}
By the minimization of the latter estimate in $\delta $ we obtain that
\begin{equation*}
\left\| Q_{N}\right\| _{2}\leq \sqrt{3}\pi ^{\frac{\nu }{2+\nu }}A^{\frac{1}{%
2+\nu }}N^{\frac{2}{2+\nu }}\text{ }.
\end{equation*}
\endproof%

\section{Asymptotic form of the maximal eigenvalue of a truncated
correlation matrix}

Let us denote by $\lambda _{m}\left( N\right) $ the maximal eigenvalue of a
positive (i.e., positive definite) matrix $Q_{N}$, $\lambda _{m}\left(
N\right) =\left\| Q_{N}\right\| $. The condition of Theorem \ref{t1} admits
the following weakening.

\begin{theorem}
\label{t2} The process $\xi \left( t\right) $ is stable if and only if
\begin{equation*}
\underset{N\rightarrow \infty }{\lim }\frac{\lambda _{m}\left( N\right) }{N}%
~\left( =\underset{N\rightarrow \infty }{\lim }\frac{1}{N}\left\|
Q_{N}\right\| \right) =0\text{ }.
\end{equation*}
\end{theorem}

\proof%
The Hilbert norm $\left\| A\right\| $ of any square matrix (or any nuclear
operator $A)$ satisfies the inequality
\begin{equation}
\left\| A\right\| \geq \frac{\left\| A\right\| _{2}}{\left\| A\right\| _{1}}%
\text{ },  \label{d0}
\end{equation}
where $\left\| A\right\| _{1}$ and $\left\| A\right\| _{2}$ are the nuclear
and Hilbert-Schmidt norms, respectively. Recall that $\left\| A\right\| $
for positive definite $A$ coincides with the maximal eigenvalue $\lambda
_{m}\left( A\right) $ of $A$ and for $A\geq 0$ (\ref{d0}) takes the form
\begin{equation}
\lambda _{m}\left( A\right) \geq \frac{\mathrm{Tr}A^{2}}{\mathrm{Tr}A}\text{
}.  \label{d1}
\end{equation}
The application of (\ref{d1}) and (\ref{a7}) to $Q_{N}$ by virtue of the
equality $\mathrm{Tr}Q_{N}=Nb\left( 0\right) =N$ gives
\begin{equation}
\underset{N\rightarrow \infty }{\underline{\lim }}\frac{\lambda _{m}\left(
N\right) }{N}\geq \frac{1}{2}\sum\limits_{\alpha }m_{\alpha }^{2}\text{ }.
\label{d2}
\end{equation}
On the other hand,
\begin{equation*}
\lambda _{m}\left( N\right) \leq \sqrt{\mathrm{Tr}Q_{N}^{2}}\text{ }.
\end{equation*}
Hence,
\begin{equation}
\underset{N\rightarrow \infty }{\overline{\lim }}\frac{\lambda _{m}\left(
N\right) }{N}\leq \sqrt{\frac{1}{2}\sum\limits_{\alpha }m_{\alpha }^{2}}%
\text{ }.  \label{d3}
\end{equation}
\endproof%

The inequalities (\ref{d2}), (\ref{d3}) can be specified.

\begin{theorem}
\label{t3} Given the sequence of maximal eigenvalues (norms) $\left\{
\lambda _{m}\left( N\right) \right\} $ of positive definite T\"{o}plitz
matrices $\left\{ Q_{N}\right\} $, generated by a non-negative measure $%
d\sigma \left( \theta \right) $ (\ref{aa2}), it holds that
\begin{equation*}
\underset{N\rightarrow \infty }{\lim }\frac{\lambda _{m}\left( N\right) }{N}=%
\underset{\alpha }{\max }m_{\alpha }\text{ .}
\end{equation*}
\end{theorem}

\proof%
The T\"{o}plitz matrix $Q_{N}$ generated by the non-decreasing function (\ref
{aa2}) is the sum of non-negative T\"{o}plitz matrices $Q_{N}^{\left(
c\right) }$ and $Q_{N}^{\left( d\right) }$, generated by non-decreasing
functions $\sigma _{c}$ and $\sigma _{d}$, respectively. Let us denote by $%
\lambda _{m}^{\left( c\right) }\left( N\right) $ and $\lambda _{m}^{\left(
d\right) }\left( N\right) $ the maximal eigenvalues (norms) of the matrices $%
Q_{N}^{\left( c\right) }$ and $Q_{N}^{\left( d\right) }$, respectively.
Since $Q_{N}$ $\geq Q_{N}^{\left( d\right) }$, then
\begin{equation*}
\lambda _{m}^{\left( d\right) }\left( N\right) \leq \lambda _{m}\left(
N\right) =\left\| Q_{N}^{\left( d\right) }+Q_{N}^{\left( c\right) }\right\|
\leq \left\| Q_{N}^{\left( d\right) }\right\| +\left\| Q_{N}^{\left(
c\right) }\right\| =\lambda _{m}^{\left( d\right) }\left( N\right) +\lambda
_{m}^{\left( c\right) }\left( N\right) \text{ }.
\end{equation*}
By virtue of Theorem \ref{t2}, $\lambda _{m}^{\left( c\right) }\left(
N\right) =o\left( N\right) $. Hence it remains to prove that
\begin{equation}
\underset{N\rightarrow \infty }{\lim }\frac{\lambda _{m}^{\left( d\right)
}\left( N\right) }{N}=\underset{\alpha }{\max }m_{\alpha }\text{ }.
\label{d5}
\end{equation}
To this end, let us consider first the case of $\sigma _{d}\left( \theta
\right) $ having only a finite number $2s$ of points of growth. We will not
use in this proof the fact that the jump points of $\sigma _{d}\left( \theta
\right) $ are located symmetrically with respect to the point $\theta =0$.
The T\"{o}plitz matrix $Q_{N}^{\left( d\right) }=\left( b_{d}\left(
j-k\right) \right) _{j,k=0}^{N-1},$ $2s\ll N$, generated by $\sigma _{d}$,
can be represented in this case in the form
\begin{equation}
Q_{N}^{\left( d\right) }=\sum\limits_{\alpha =1}^{2s}m_{\alpha }\left( \cdot
,\mathbf{e}_{\alpha }\right) \mathbf{e}_{\alpha }\text{ },  \label{a10}
\end{equation}
where
\begin{equation*}
\left( \cdot ,\mathbf{e}_{\alpha }\right) \mathbf{e}_{\alpha }=\left(
\mathbf{\exp }\left( i\left( j-k\right) \theta _{\alpha }\right) \right)
_{j,k=0}^{N-1}
\end{equation*}
are $N\times N$ matrices of unit rank, so that $Q_{N}$ transforms a $N\times
1$ column vector $\mathbf{x}=\left( x_{j}\right) _{j=0}^{N-1}$ into
\begin{equation}
Q_{N}^{\left( d\right) }\mathbf{x}=\sum\limits_{\alpha =1}^{2s}m_{\alpha
}\left( \mathbf{x},\mathbf{e}_{\alpha }\right) \mathbf{e}_{\alpha }\text{ },
\label{a11}
\end{equation}
where $\left( \mathbf{\cdot },\mathbf{\cdot }\right) $ is the scalar product
in the linear space of $N\times 1$ column vectors $\Bbb{{C}_{N}}$ defined in
a standard way:
\begin{equation*}
\left( \mathbf{x},\mathbf{y}\right) =\sum\limits_{j=0}^{N-1}x_{j}\overline{y}%
_{j},\;\mathbf{x}=\left( x_{j}\right) _{j=0}^{N-1},\,\mathbf{y}=\left(
y_{j}\right) _{j=0}^{N-1}\text{ }.
\end{equation*}

Notice that the vectors $\left\{ \mathbf{e}_{\alpha }\right\} $ are linearly
independent. Indeed, suppose that there is a set of complex numbers $\left\{
z_{\alpha }\right\} $ such that
\begin{equation}
\sum\limits_{\alpha =1}^{s}z_{\alpha }\mathbf{e}_{\alpha }=0\text{ }.
\label{a12}
\end{equation}
Due to (\ref{a12}), the numbers $z_{\alpha }$ satisfy the homogeneous system
\begin{equation*}
\sum\limits_{\alpha =1}^{s}\exp \left( ik\theta _{\alpha }\right) z_{\alpha
}=0,\,\;k=0,1,...,2s-1.
\end{equation*}
But the determinant of this system is the Van der Monde determinant, which
vanishes if and only if among the numbers $\left\{ \exp \left( ik\theta
_{\alpha }\right) \right\} $ there are equals. The latter is impossible by
our assumption. Hence all $z_{\alpha }=0$.

Let $\lambda $ be a non-zero eigenvalue of $Q_{N}^{\left( d\right) }$ and $%
\mathbf{h}_{\lambda }$ be a corresponding non-zero eigenvector:
\begin{equation}
\sum\limits_{\alpha =1}^{s}m_{\alpha }\left( \mathbf{h}_{\lambda },\mathbf{e}%
_{\alpha }\right) \mathbf{e}_{\alpha }=\lambda \mathbf{h}_{\lambda }\text{ }.
\label{b2}
\end{equation}
By (\ref{b2}) $\mathbf{h}_{\lambda }$ admits the representation:
\begin{equation*}
\mathbf{h}_{\lambda }=\sum\limits_{\alpha =1}^{s}z_{\alpha }\mathbf{e}%
_{\alpha }\text{ },
\end{equation*}
where $z_{\alpha }$ are some complex numbers, not all of which are equal to
zero. Put
\begin{equation*}
\eta _{\alpha }=\sqrt{m_{\alpha }}\left( \mathbf{h}_{\lambda },\mathbf{e}%
_{\alpha }\right) \text{ }.
\end{equation*}
By virtue of (\ref{b2}), not all numbers $\eta _{\alpha }=0$. Taking the
scalar products of both sides of (\ref{b2}) with all vectors $\sqrt{%
m_{\alpha }}\mathbf{e}_{\alpha }$, we obtain the following homogeneous
system for $\eta _{\alpha }$:
\begin{equation}
\sum\limits_{\alpha ^{\prime }=1}^{s}\sqrt{m_{\alpha }m_{\alpha ^{\prime }}}%
\left( \mathbf{e}_{\alpha ^{\prime }},\mathbf{e}_{\alpha }\right) \eta
_{\alpha ^{\prime }}=\lambda \eta _{\alpha }\text{ }.  \label{b3}
\end{equation}
Thus, the non-zero eigenvalues of $Q_{N}^{\left( d\right) }$ coincide, with
account of their multiplicities, with the eigenvalues of the $2s\times 2s$
Hermitian positive definite matrix
\begin{equation}
A_{N}=\left( \sqrt{m_{\alpha }m_{\alpha ^{\prime }}}\left( \mathbf{e}%
_{\alpha ^{\prime }},\mathbf{e}_{\alpha }\right) \right) _{\alpha ,\alpha
^{\prime }=1}^{2s}\text{ }.  \label{d6}
\end{equation}
Notice that by definition of the vectors $\mathbf{e}_{\alpha }$, we have
\begin{equation}
\left( \mathbf{e}_{\alpha ^{\prime }},\mathbf{e}_{\alpha }\right) =\frac{%
\exp \left( iN\left( \theta _{\alpha ^{\prime }}-\theta _{\alpha }\right)
\right) -1}{\exp \left( i\left( \theta _{\alpha ^{\prime }}-\theta _{\alpha
}\right) \right) -1},\;\alpha ^{\prime }\neq \alpha ;\;\left( \mathbf{e}%
_{\alpha },\mathbf{e}_{\alpha }\right) =N\text{ }.  \label{d7a}
\end{equation}
Hence, the matrix $A_{N}$ is the sum of the diagonal matrix
\begin{equation*}
A_{1,N}:=\left( Nm_{\alpha }\delta _{\alpha \alpha ^{\prime }}\right)
_{\alpha ,\alpha ^{\prime }=1}^{2s}
\end{equation*}
and the Hermitian matrix $A_{2,N}$ with zero diagonal elements and
non-diagonal elements $\sqrt{m_{\alpha }m_{\alpha ^{\prime }}}\left( \mathbf{%
e}_{\alpha ^{\prime }},\mathbf{e}_{,\alpha }\right) ,$ $\alpha \neq \alpha
^{\prime }$. By (\ref{d7a}) the non-diagonal elements of $A_{2,N}$ are
uniformly bounded:
\begin{equation*}
\left| \sqrt{m_{\alpha }m_{\alpha ^{\prime }}}\left( \mathbf{e}_{\alpha
^{\prime }},\mathbf{e}_{\alpha }\right) \right| \leq 2\left( \underset{%
\alpha ^{\prime }\neq \alpha }{\max }\left| \theta _{\alpha ^{\prime
}}-\theta _{\alpha }\right| ^{-1}\right) \left( \underset{\alpha }{\max }%
m_{\alpha }\right) \text{ },
\end{equation*}
and, hence,
\begin{equation*}
\left\| A_{2,N}\right\| \leq \left( 4s-2\right) \left( \underset{\alpha
^{\prime }\neq \alpha }{\max }\left| \theta _{\alpha ^{\prime }}-\theta
_{\alpha }\right| ^{-1}\right) \left( \underset{\alpha }{\max }m_{\alpha
}\right) \text{ }.
\end{equation*}
Therefore,
\begin{equation*}
\left[ N-\left( 4s-2\right) \left( \underset{\alpha ^{\prime }\neq \alpha }{%
\max }\left| \theta _{\alpha ^{\prime }}-\theta _{\alpha }\right|
^{-1}\right) \right] \left( \underset{\alpha }{\max }m_{\alpha }\right) \leq
\left\| A_{1,N}\right\| -\left\| A_{2,N}\right\| \leq
\end{equation*}
\begin{equation*}
\left\| A_{N}\right\| \leq \left\| A_{1,N}\right\| +\left\| A_{2,N}\right\|
\leq \left[ N+\left( 4s-2\right) \left( \underset{\alpha ^{\prime }\neq
\alpha }{\max }\left| \theta _{\alpha ^{\prime }}-\theta _{\alpha }\right|
^{-1}\right) \right] \left( \underset{\alpha }{\max }m_{\alpha }\right)
\text{ }.
\end{equation*}
We see that
\begin{equation}
\lambda _{m}^{\left( d\right) }\left( N\right) =\left\| A_{N}\right\|
\underset{N\rightarrow \infty }{=}N\cdot \underset{\alpha }{\max }m_{\alpha
}+O\left( 1\right) \text{ }.  \label{d7b}
\end{equation}

To prove the relation (\ref{d5}) for a non-decreasing step function $\sigma
_{d}\left( \theta \right) ,\;\sigma _{d}\left( \pi \right) -\sigma
_{d}\left( -\pi \right) \leq 1$, having infinitely many points of jump, we
take a small $\varepsilon >0$ and split $\sigma _{d}\left( \theta \right) $
into a sum $\sigma _{1d}\left( \theta \right) +\sigma _{2d}\left( \theta
\right) $ of non-decreasing step functions $\sigma _{1,d}\left( \theta
\right) $ and $\sigma _{2,d}\left( \theta \right) $, where, as before, $%
\sigma _{1d}\left( \theta \right) $ has a finite number of jump points and $%
\sigma _{2d}\left( \theta \right) $ is such that
\begin{equation*}
\int\limits_{-\pi }^{\pi }d\sigma _{2,d}\left( \theta \right) <\varepsilon <%
\underset{\alpha }{\max }m_{\alpha }\text{ }.
\end{equation*}
With respect to this split, we represent the T\"{o}plitz matrix $%
Q_{N}^{\left( d\right) }$ as the sum $Q_{N}^{\left( 1,d\right)
}+Q_{N}^{\left( 2,d\right) }$ of non-negative T\"{o}plitz matrices generated
by $\sigma _{1,d}\left( \theta \right) $ and $\sigma _{2,d}\left( \theta
\right) $, respectively. Notice that by construction
\begin{equation}
\left\| Q_{N}^{\left( 2,d\right) }\right\| \leq \mathrm{Tr}Q_{N}^{\left(
2,d\right) }<N\varepsilon \text{ }.  \label{d7c}
\end{equation}
Besides,
\begin{equation*}
\lambda _{m}^{\left( 1,d\right) }\left( N\right) =\left\| Q_{N}^{\left(
1,d\right) }\right\| \leq \lambda _{m}^{\left( d\right) }\left( N\right)
\leq \lambda _{m}^{\left( 1,d\right) }\left( N\right) +\left\| Q_{N}^{\left(
2,d\right) }\right\| \text{ }.
\end{equation*}
Applying the estimate (\ref{d7b}) to $Q_{N}^{\left( 1,d\right) }$ and taking
into account the inequality (\ref{d7c}) for $N\rightarrow \infty $ yields
\begin{equation*}
N\cdot \underset{\alpha }{\max }m_{\alpha }+O\left( 1\right) =\lambda
_{m}^{\left( 1,d\right) }\left( N\right) \leq \lambda _{m}^{\left( d\right)
}\left( N\right) \leq N\cdot \left( \underset{\alpha }{\max }m_{\alpha
}+\varepsilon \right) +O\left( 1\right) \text{ }.
\end{equation*}
Finally,
\begin{equation*}
\underset{\alpha }{\max }m_{\alpha }\leq \underset{N\rightarrow \infty }{%
\underline{\lim }}\frac{\lambda _{m}^{\left( d\right) }\left( N\right) }{N}%
\leq \underset{N\rightarrow \infty }{\overline{\lim }}\frac{\lambda
_{m}^{\left( d\right) }\left( N\right) }{N}\leq \underset{\alpha }{\max }%
m_{\alpha }+\varepsilon \text{ },
\end{equation*}
where $\varepsilon >0$ can be taken arbitrarily small.%
\endproof%

\begin{remark}
\label{t4} For the Hilbert norm $\left\| C\right\| $ of any $N\times N$
matrix $C=\left( c_{jk}\right) _{j,k=0}^{N-1}$ the following estimate holds:
\begin{equation}
\left\| C\right\| \leq \max \left\{ \underset{j}{\max }\sum%
\limits_{p=0}^{N-1}\left| c_{jp}\right| ,\,\underset{k}{\max }%
\sum\limits_{p=0}^{N-1}\left| c_{pk}\right| \right\} .  \label{d8}
\end{equation}
If $C=\left( c_{j-k}\right) _{j,k=0}^{N-1}$ is a T\"{o}plitz matrix, then by
(\ref{d8})
\begin{equation}
\left\| C\right\| \leq \left| c_{0}\right| +\sum\limits_{p=1}^{N-1}\left(
\left| c_{p}\right| +\left| c_{-p}\right| \right) \text{ }.  \label{d9}
\end{equation}
Thus for the T\"{o}plitz matrices $Q_{N}$, which are under consideration
here, Theorem \ref{t3} by virtue of (\ref{d9}), gives
\begin{equation}
\frac{1}{2}\underset{\alpha }{\max }m_{\alpha }=\underset{N\rightarrow
\infty }{\lim }\frac{1}{2N}\left\| Q_{N}\right\| \leq \underset{N\rightarrow
\infty }{\underline{\lim }}\frac{1}{N}\sum\limits_{p=0}^{N-1}\left| b\left(
p\right) \right| \text{ }.  \label{d10}
\end{equation}
Therefore the condition
\begin{equation}
\underset{N\rightarrow \infty }{\lim }\frac{1}{N}\sum\limits_{p=0}^{N-1}%
\left| b\left( p\right) \right| =0  \label{d10a}
\end{equation}
is sufficient for the process stability.
\end{remark}

\section{Extension to continuous time processes}

Real signals are, certainly, continuous time processes, $\xi \left( t\right)
$. The correlation function $b\left( t\right) $ of a process having a finite
second moment $\left\langle \xi ^{2}\left( t\right) \right\rangle $ is a
Hermitian positive function. As such, it admits the representation
\begin{equation}
b\left( t\right) =\int_{-\infty }^{\infty }\exp \left( i\lambda t\right)
d\vartheta \left( \lambda \right) \text{ },  \label{c1}
\end{equation}
where $\vartheta \left( \lambda \right) $ is a bounded non-decreasing
function on the real axis. Like for the discrete time processes, $\vartheta
\left( \lambda \right) $ can be represented, in general, as the sum
\begin{equation*}
\vartheta \left( \lambda \right) =\vartheta _{c}\left( \lambda \right)
+\vartheta _{d}\left( \lambda \right)
\end{equation*}
of a non-decreasing continuous function $\vartheta _{c}\left( \lambda
\right) $ and a non-decreasing step function $\vartheta _{d}\left( \lambda
\right) $, and we call the process \textit{stable} if $\vartheta \left(
\lambda \right) $ is continuous and \textit{unstable} otherwise. To
investigate the instability characteristics of a continuous time process, we
consider instead of the T\"{o}plitz matrices $Q_{N}$, the set of
non-negative integral operators
\begin{equation}
\left( B_{T}f\right) \left( t\right) =\int_{0}^{T}b\left( t-s\right) f\left(
s\right) ds,\;0<T<\infty \text{ },  \label{d12}
\end{equation}
in the Hilbert spaces $L_{2}\left( 0,T\right) $. Since $b\left( t\right) $
is a continuous function, all these operators are nuclear and their nuclear
and Hilbert-Schmidt norms $\left\| B_{T}\right\| _{1}$ and $\left\|
B_{T}\right\| _{2}$ are given by the expressions
\begin{equation}
\left\| B_{T}\right\| _{1}=Tb\left( 0\right) =\int_{-\infty }^{\infty
}d\vartheta \left( \lambda \right) \text{ },  \label{e1}
\end{equation}
\begin{equation*}
\left\| B_{T}\right\| _{2}=\sqrt{2T\int\limits_{0}^{T}\left( 1-\frac{t}{T}%
\right) \left| b\left( t\right) \right| ^{2}dt}=
\end{equation*}
\begin{equation}
\sqrt{\int\limits_{-\infty }^{\infty }\int\limits_{-\infty }^{\infty }\frac{4%
}{\left( \lambda -\lambda ^{\prime }\right) ^{2}}\sin ^{2}\frac{\left(
\lambda -\lambda ^{\prime }\right) T}{2}d\vartheta \left( \lambda ^{\prime
}\right) d\vartheta \left( \lambda \right) }\text{ }.  \label{e2}
\end{equation}
Let us denote, as before, by $\left\{ m_{\alpha }\right\} $ the set of jumps
of $\vartheta \left( \lambda \right) $. Using (\ref{e1}), (\ref{e2}) and the
arguments similar to those employed in the proofs of Theorems \ref{t1} and
\ref{t3}, we obtain the following criteria of stability of a continuous time
process $\xi \left( t\right) $.

\begin{theorem}
\label{t5} A stationary continuous time process $\xi \left( t\right) $ is
stable if and only if its correlation function $b\left( t\right) $ possesses
the property:
\begin{equation*}
\underset{T\rightarrow \infty }{\lim }\frac{2}{T}\int\limits_{0}^{T}\left( 1-%
\frac{t}{T}\right) \left| b\left( t\right) \right| ^{2}dt=0\text{ }.
\end{equation*}
Otherwise,
\begin{equation*}
\underset{T\rightarrow \infty }{\lim }\frac{2}{T}\int\limits_{0}^{T}\left( 1-%
\frac{t}{T}\right) \left| b\left( t\right) \right|
^{2}dt=\sum\limits_{\alpha }m_{\alpha }^{2}.
\end{equation*}
\end{theorem}

\begin{theorem}
\label{t6} A stationary continuous time process $\xi \left( t\right) $ is
stable if and only if the Hilbert norms $\left\| B_{T}\right\| $ of integral
operators (\ref{d12}), where $b\left( t\right) $ is the correlation function
of the process, are such that
\begin{equation*}
\underset{T\rightarrow \infty }{\lim }\frac{1}{T}\left\| B_{T}\right\| =0%
\text{ }.
\end{equation*}
Otherwise,
\begin{equation}
\underset{T\rightarrow \infty }{\lim }\frac{1}{T}\left\| B_{T}\right\| =%
\underset{\alpha }{\max }m_{\alpha }\text{ }.  \label{e3}
\end{equation}
\end{theorem}

\begin{remark}
\label{t8} For the norm of the integral operator $B_{T}$ the following
estimate:
\begin{equation*}
\left\| B_{T}\right\| \leq 2\int\limits_{0}^{T}\left| b\left( t\right)
\right| dt
\end{equation*}
holds. As it stems from (\ref{e3}) the relation
\begin{equation*}
\underset{T\rightarrow \infty }{\lim }\frac{1}{T}\int\limits_{0}^{T}\left|
b\left( t\right) \right| dt=0\text{ },
\end{equation*}
guarantees the stability of the process $\xi \left( t\right) $ $.$
\end{remark}

\section{Application to processing of random signals}

\subsection{Detection of quasi-periodic components in a random signal}

The jumps of the spectral distribution function $\vartheta \left( \lambda
\right) $ of a stationary process is, in general, a sign of appearance of
undamped oscillation components in a signal and the points of discontinuity
of $\vartheta \left( \lambda \right) $ are either the frequencies of such
components themselves or directly related to them. Notice that, due to
physical reasons, the measurement of $\xi \left( t\right) $ is possible only
at discrete moments of time with a step $\Delta $. If $t$ in (\ref{c1}) is
an integer multiple of $\Delta $, then it is evident that
\begin{eqnarray}
b\left( t\right) &=&\int_{-\Omega }^{\Omega }\exp \left( it\theta \right)
d\sigma \left( \theta \right) ,\;\Omega =\frac{\pi }{\Delta }\text{ },
\label{c2} \\
\sigma \left( \theta \right) &=&\sum\limits_{n=-\infty }^{\infty }\left[
\vartheta \left( \theta +2n\Omega \right) -\vartheta \left( 2n\Omega -\Omega
\right) \right] ,\;-\Omega \leq \theta <\Omega \text{ }.  \notag
\end{eqnarray}
The function $\sigma \left( \theta \right) $ is bounded and non-decreasing
in the interval $[-\Omega ,\Omega ]$. If $\vartheta \left( \lambda \right) $
loses its continuity at the points $\pm \lambda _{1},\pm \lambda _{2},$ $.$ $%
.$ $.$ , then $\sigma \left( \theta \right) $ has a non-void set of
discontinuity points
\begin{equation}
\left\{ \pm \theta _{j}^{\prime }=\left( \mathbf{E}\left( \frac{\pm \lambda
_{j}}{2\Omega }+\frac{1}{2}\right) -\frac{1}{2}\right) 2\Omega \right\}
\subset [-\Omega ,\Omega ]\text{ },  \label{c3}
\end{equation}
where $\mathbf{E}\left( x\right) $ is the fractional part of the number $x$.
(In general, $\pm \theta _{j_{1}}^{\prime }$ coincides with every $\pm
\theta _{j_{2}}^{\prime }$ such that $\lambda _{j_{1}}-\lambda _{j_{2}}$ is
a multiple of $2\Omega $.) Therefore, in general, the jump of $\sigma \left(
\theta \right) $ at a point $\theta _{j}^{\prime }$ is the sum of the jumps
of $\vartheta \left( \lambda \right) $ at all co-images of $\theta
_{j}^{\prime }$ under the mapping (\ref{c3}).) Taking $\Delta $ as the time
measurement unit, we return to the representation (\ref{a1}) of $b\left(
t\right) $ for integer $t$. Thus, the spectral distribution function $\sigma
\left( \theta \right) $ inherits all discontinuities of $\vartheta \left(
\lambda \right) $ from the interval $[-\Omega ,\Omega ]$ and also may get
new ones at the points (calculated according to (\ref{c3})) related to the
discontinuity points of $\vartheta \left( \lambda \right) $ outside this
interval. We see that the spectral distribution function for the discrete
time process obtained in such a way from a continuous time process, has a
non-trivial component $\sigma _{d}$ if and only if the corresponding
spectral distribution function of the initial discrete time process has
non-zero jumps on some set of points. In other words, the values of a random
continuous time process measured at discrete moments of time, form a stable
discrete time process if and only if the initial process is stable.

The correlation function of the discrete time process delivers not only the
described gauge of instability of the process, but also the following tool
for the detection of the points $\left\{ \pm \theta _{\alpha }\right\} $,
which are the discontinuity points of $\sigma \left( \theta \right) $. Put
\begin{equation*}
\Theta _{N}\left( \theta \right) =b\left( 0\right)
+2\sum\limits_{k=0}^{N-1}\left( 1-\frac{k}{N}\right) b\left( k\right) \cos
k\theta =
\end{equation*}
\begin{equation}
=\int\limits_{-\pi }^{\pi }\frac{\sin ^{2}\frac{1}{2}N\left( \theta -\theta
^{\prime }\right) }{N^{2}\sin ^{2}\frac{1}{2}\left( \theta -\theta ^{\prime
}\right) }d\sigma \left( \theta ^{\prime }\right) \text{ }.  \label{c4}
\end{equation}
It is not difficult to see that
\begin{equation}
\underset{N\rightarrow \infty }{\lim }\frac{1}{N}\Theta _{N}\left( \theta
\right) =\sigma \left( \theta +0\right) -\sigma \left( \theta -0\right)
\text{ }.  \label{c5}
\end{equation}
Further, take a sufficiently large $N\gg 1$ and split the interval $[-\pi
,\pi ]$ into equal segments of longitude $\delta $ such that $N\delta
\lesssim 1$. Let $\sigma \left( \theta \right) $ have a jump $m_{\alpha }$
within the interval $\left( l\delta ,\left( l+1\right) \delta \right) $.
Since
\begin{equation*}
\frac{5}{6}\left| x\right| \leq \left| x\right| \left( 1-\frac{x^{2}}{6}%
\right) \leq \left| \sin x\right| \leq \left| x\right| \text{ },
\end{equation*}
then, for $\left| \theta -l\delta \right| \sim \delta $, we have
\begin{equation*}
\Theta _{N}\left( \theta \right) \geq \int\limits_{l\delta }^{\left(
l+1\right) \delta }\frac{\sin ^{2}\frac{1}{2}N\left( \theta -\theta ^{\prime
}\right) }{N^{2}\sin ^{2}\frac{1}{2}\left( \theta -\theta ^{\prime }\right) }%
d\sigma \left( \theta ^{\prime }\right) \geq \frac{5}{6}N\int\limits_{l%
\delta }^{\left( l+1\right) \delta }d\sigma \left( \theta ^{\prime }\right)
\geq \frac{5}{6}Nm_{\alpha }\text{ }.
\end{equation*}

On the other hand, if the continuous part $\sigma _{c}\left( \theta \right) $
of $\sigma \left( \theta \right) $ satisfies the condition of Proposition
\ref{pr}, then one can use the arguments similar to those employed in the
proof of this proposition to show that
\begin{equation*}
\Theta _{N}\left( \theta \right) =O\left( N^{\frac{2}{2+\nu }}\right)
\end{equation*}
at the points remote from the jumps of $\sigma \left( \theta \right) $.
Hence, the values of $\Theta _{N}\left( l\delta \right) $, where $l$ is an
integer which satisfies the condition $-\pi \leq l\delta \leq \pi $, at the
distances $\sim \delta $ from the jump points of $\sigma \left( \theta
\right) $, must be visible as larger than those at the distances$\sim $ $%
l_{0}\delta $, where $Nl_{0}\delta \sim 1$.

The assertion of Theorems \ref{t1} and \ref{t3} can be used for the
detection of symptoms of emerging instabilities of a random process, which
can be considered as stationary for long time intervals. The method consists
in the construction of the correlation function of the process from a piece
of its time series from the beginning of observation to a rather far off
moment of time $\Upsilon $ in the future. Set, as usually,
\begin{equation*}
b\left( k\right) =\frac{1}{\Upsilon -k}\sum\limits_{p=0}^{\Upsilon -k}\xi
\left( p+k\right) \xi \left( p\right) -m^{2},\;m=\frac{1}{\Upsilon }%
\sum\limits_{p=0}^{\Upsilon }\xi \left( p\right) \text{ },
\end{equation*}
and compute, for a sufficiently large $N\ll \Upsilon $, the numbers
\begin{equation}
\frac{1}{N^{2}}\left\| Q_{N}\right\| _{2}^{2}=\frac{1}{N}+\frac{2}{N}%
\sum\limits_{k=1}^{N-1}\left( 1-\frac{k}{N}\right) \frac{b^{2}\left(
k\right) }{b^{2}\left( 0\right) },\;\frac{1}{N}\sum\limits_{p=0}^{N-1}\left|
\frac{b\left( p\right) }{b\left( 0\right) }\right| \text{ },  \label{a8}
\end{equation}
or the numbers
\begin{equation*}
\frac{1}{T}\int\limits_{0}^{T}\left( 1-\frac{t}{T}\right) \left| b\left(
t\right) \right| ^{2}dt,\;\frac{1}{T}\int\limits_{0}^{T}\left| b\left(
t\right) \right| dt
\end{equation*}
for a continuous time process. An explicit tendency of any of these numbers
to be bounded, for increasing $N$, from below by certain positive numbers,
is a serious evidence of the process instability.

The following example demonstrates that the manifestation of such a tendency
begins the sooner in $N$ the larger is the contribution of the oscillating
components generated by $d\sigma _{d}\left( \theta \right) $ into $b\left(
0\right) $.

Let the correlation function of a stationary random process be given by the
expressions:
\begin{eqnarray}
b\left( 0\right) &=&1;\;b\left( k\right) =b\left( -k\right)
=\sum\limits_{\alpha =1}^{s}m_{\alpha }\,\cos k\theta _{\alpha },\;k=1,2,...%
\text{ },  \label{a9} \\
1 &\leq &s<\infty ,\;0<\sum\limits_{\alpha =1}^{s}m_{\alpha }<1,\;0<\theta
_{1},...,\theta _{s}<\pi \text{ }.  \notag
\end{eqnarray}
The spectral distribution function $\sigma \left( \theta \right) $ of such a
process is the sum of the spectral distribution function $\sigma _{c}\left(
\theta \right) $ of the ''white noise'',
\begin{equation*}
d\sigma _{c}\left( \theta \right) =\frac{p}{2\pi }d\theta
,\;p=1-\sum\limits_{\alpha =1}^{s}m_{\alpha }\text{ },
\end{equation*}
and the step function $\sigma _{d}\left( \theta \right) $, the jump points
of which are $\left\{ \pm \theta _{\alpha }\right\} $, and
\begin{equation*}
\sigma _{d}\left( -\theta _{\alpha }+0\right) -\sigma _{d}\left( -\theta
_{\alpha }-0\right) =\sigma _{d}\left( \theta _{\alpha }+0\right) -\sigma
_{d}\left( \theta _{\alpha }-0\right) =\frac{1}{2}m_{\alpha }\text{ }.
\end{equation*}
To calculate the right-hand side of (\ref{a4}) in this special case, notice
that
\begin{equation*}
\frac{1}{2\pi }\int\limits_{-\pi }^{\pi }\frac{\sin ^{2}\frac{1}{2}N\left(
\theta -\theta ^{\prime }\right) }{\sin ^{2}\frac{1}{2}\left( \theta -\theta
^{\prime }\right) }d\theta ^{\prime }\equiv N\text{ }.
\end{equation*}
Therefore we see that now $\Phi _{N}\left( \theta \right) =pN^{-1}$ and thus
\begin{equation}
\frac{1}{N^{2}}\left\| Q_{N}\right\| _{2}^{2}=\frac{1}{2}\sum\limits_{\alpha
=1}^{s}m_{\alpha }^{2}+\frac{p\left( 2-p\right) }{N}+  \label{a9a}
\end{equation}
\begin{equation*}
\frac{1}{2N^{2}}\sum\limits_{\alpha ^{\prime }\neq \alpha }m_{\alpha
}m_{\alpha ^{\prime }}\left( \frac{\sin ^{2}\frac{1}{2}N\left( \theta
_{\alpha }-\theta _{\alpha ^{\prime }}\right) }{\sin ^{2}\frac{1}{2}\left(
\theta _{\alpha }-\theta _{\alpha ^{\prime }}\right) }+\frac{\sin ^{2}\frac{1%
}{2}N\left( \theta _{\alpha }+\theta _{\alpha ^{\prime }}\right) }{\sin ^{2}%
\frac{1}{2}\left( \theta _{\alpha }+\theta _{\alpha ^{\prime }}\right) }%
\right) \text{ }.
\end{equation*}

Hence, the first term on the right hand side of (\ref{a9a}) becomes dominant
for
\begin{equation*}
N\gtrsim 2p\left( 2-p\right) \left( \sum\limits_{\alpha =1}^{s}m_{\alpha
}^{2}\right) ^{-1}\text{ }.
\end{equation*}

\subsection{Numerical results}

Let us apply the latter results to the investigation of real signals
obtained from the Forsmark 1\&2 boiling water reactor (BWR) \cite{VGMNPLCRS}%
, and also to simulated random signals.

In Fig.1 and Fig2 we display the sequences of the values
\begin{equation}
\frac{1}{N^{2}}\left\| Q_{N}\right\| _{2}^{2}=\frac{1}{N^{2}}%
\sum\limits_{j,k=0}^{N-1}b^{2}\left( j-k\right)  \label{p1}
\end{equation}
for the T\"{o}plitz matrices constructed for two different shots of real
signals obtained in the Forsmark BWR, named in \cite{VGMNPLCRS} as \textit{%
c4\_lprm3} and \textit{c4\_lprm22}. The sampling interval of these signals
equal $0$ $.08$ $s$, and both of them consist of 4209 points. In \cite
{VGMNPLCRS} several methods were used to obtain the corresponding decay
ratio ($DR$) value. The decay ratio is a parameter related to the signal
stability: the signal is more unstable if its $DR$ is higher, see Appendix
for details. The mean values of the $DR$ given by these methods are $0$ $.90$
for \textit{c4\_lprm3} and $0$ $.51$ for \textit{c4\_lprm22}.

\begin{figure}
\begin{center}
\includegraphics[width=0.9\textwidth]{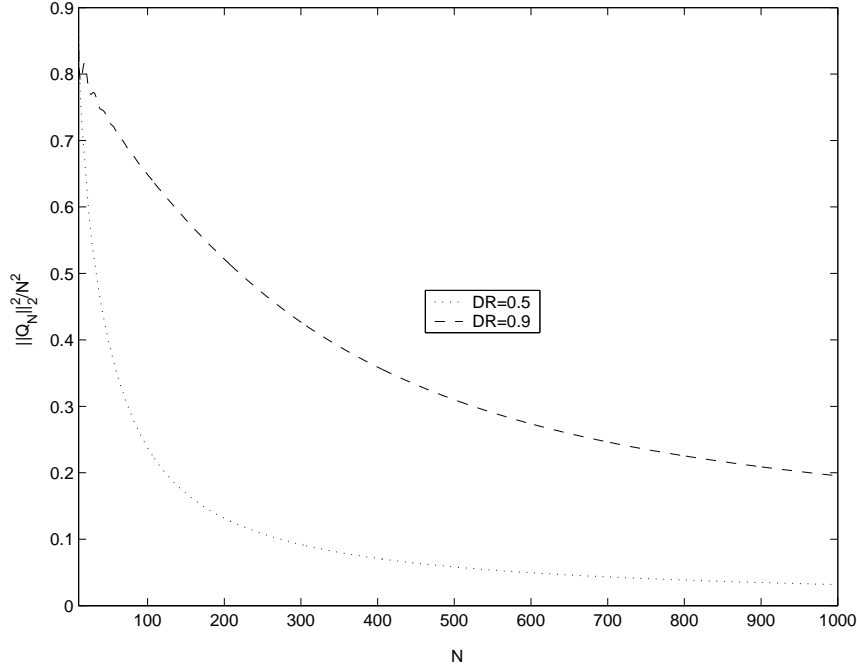}
\caption{\label{fig:1} Sequences of $\frac{1}{N^{2}}%
\left\| Q_{N}\right\| _{2}^{2}$ for the signals \textit{c4\_lprm3} and
\textit{c4\_lprm22}, linear scale of $N$.}
\end{center}
\end{figure}

\begin{figure}
\begin{center}
\includegraphics[width=0.9\textwidth]{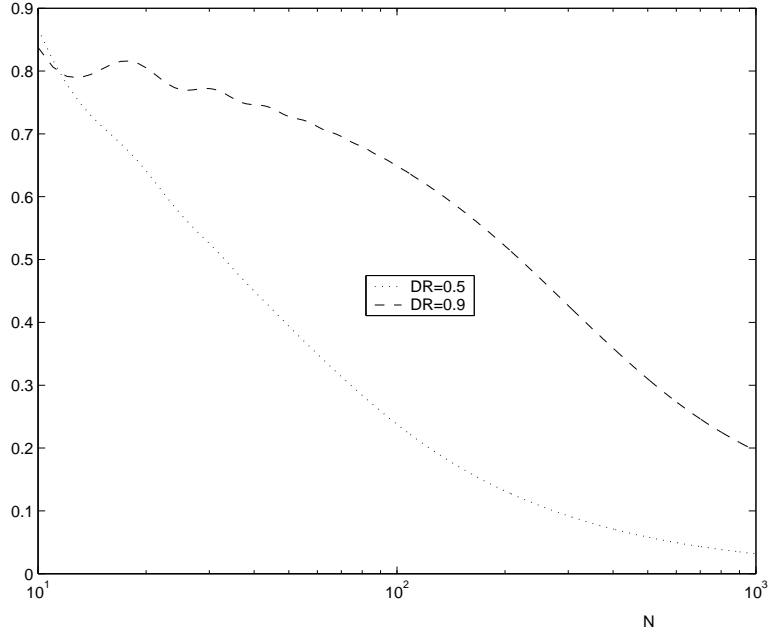}
\caption{\label{fig:2} Sequences of $\frac{1}{N^{2}}%
\left\| Q_{N}\right\| _{2}^{2}$ for the signals \textit{c4\_lprm3} and
\textit{c4\_lprm22}, logarithmic scale of $N$.}
\end{center}
\end{figure}

We observe the difference between a more unstable shot, \textit{c4\_lprm3}
and a more stable one, \textit{c4\_lprm22}. As expected, for higher $DR$,
the sequence tends to zero slower.

To determine the points of discontinuity of $\sigma \left( \theta \right) $
we have also calculated the function
\begin{equation*}
\frac{1}{N}\Theta _{N}\left( \theta \right) =\frac{b\left( 0\right) }{N}+%
\frac{2}{N}\sum\limits_{k=0}^{N-1}\left( 1-\frac{k}{N}\right) b\left(
k\right) \cos k\theta
\end{equation*}
for different values of $N$. The corresponding results are provided at Fig.3
and Fig.4, they demonstrate different behavior of $\frac{1}{N}\Theta
_{N}\left( \theta \right) $ with $N$ growing, and for different degrees of
stability measured by $DR$. In both cases the number of segments into which
the interval [-$\pi $ , $\pi $] was split, was equal to 3000. It can be
observed that for \textit{c4\_lprm22}, a more stable signal, the function $%
\frac{1}{N}\Theta _{N}\left( \theta \right) $ tends to zero with increasing $%
N$ much more rapidly.

\begin{figure}
\begin{center}
\includegraphics[width=0.9\textwidth]{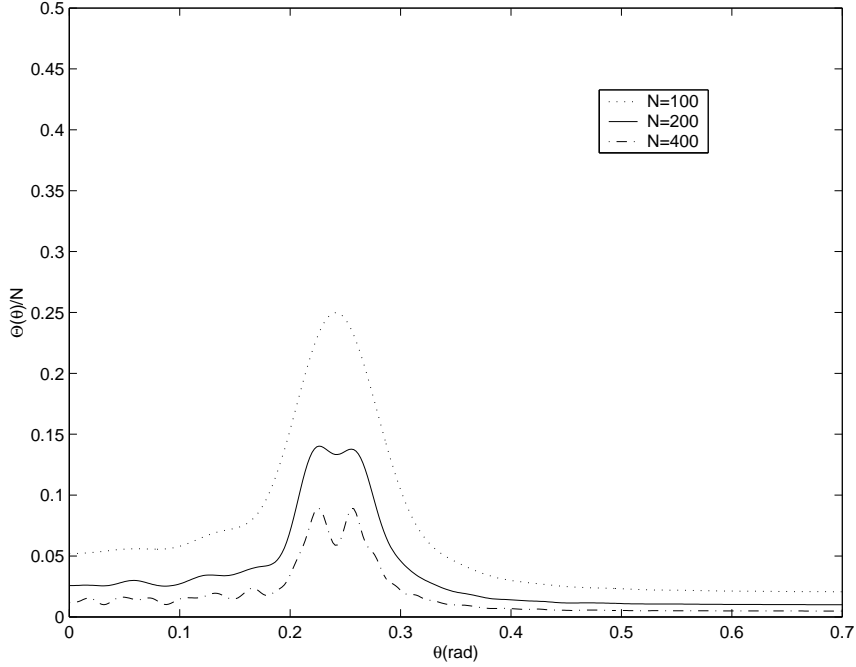}
\caption{\label{fig:3} Function $\frac{1}{N}\Theta
_{N}\left( \theta \right) $, for the signal \textit{c4\_lprm3}, $N=100,$ $%
N=300, $ $N=500$.}
\end{center}
\end{figure}

\begin{figure}
\begin{center}
\includegraphics[width=0.9\textwidth]{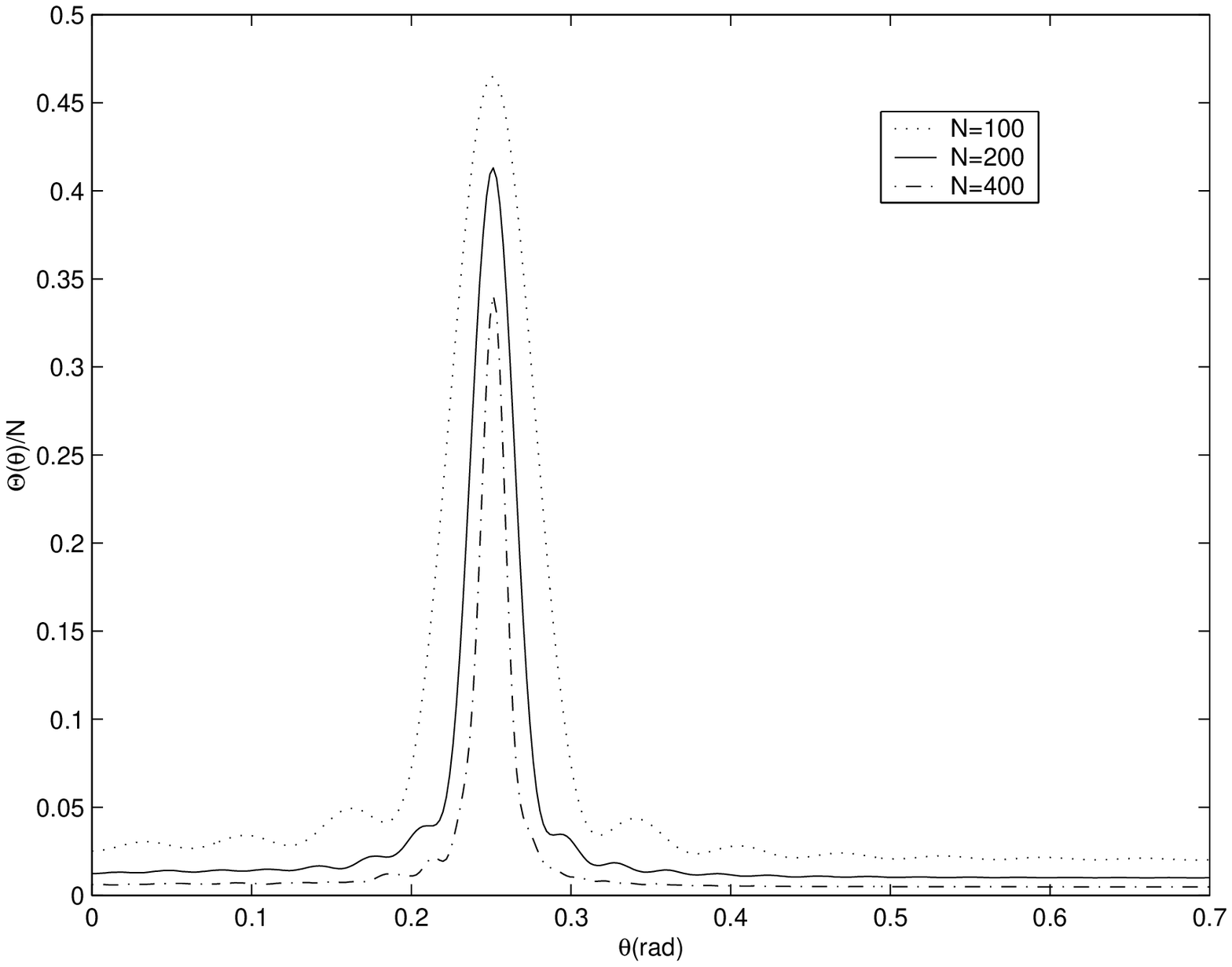}
\caption{\label{fig:4} Function $\frac{1}{N}\Theta
_{N}\left( \theta \right) $, for the signal \textit{c4\_lprm22}, $N=100,$ $%
N=300,$ $N=500$.}
\end{center}
\end{figure}

The peaks of the signal correspond to the instability frequencies. The main
peak for \textit{c4\_lprm3} is obtained for $\theta =0.24$ $rad$, which
corresponds to the frequency of $f=0.48$ $Hz$. On the other hand, for
\textit{c4\_lprm22} the main peak is located at $\theta =0.27$ $rad$ ($%
f=0.53$ $Hz$). We can compare these results with those obtained in \cite
{VGMNPLCRS}. The range of main instability frequencies obtained there is $[0.480,0.495]$
$Hz$ for \textit{c4\_lprm3} and $[0.450,0.557]$ $Hz$ for
\textit{c4\_lprm22}. The method allows to locate secondary instability
frequencies situated quite close. In Fig. 3 there is a secondary peak
located at $\theta =0.27$ $rad$, it coincides with the main peak of the
shot \textit{c4\_lprm22.} Notice also that the present approach to the
detection of instabilities is model-free.

In addition, we have also simulated signals with different degrees of
stability using the model suggested in \cite{SVM}, i.e., from the following
continuous Langevin model:
\begin{equation}
\ddot{\xi}\left( t\right) +c\dot{\xi}\left( t\right) +U\left( \xi \right)
=F\left( t\right) ,  \label{S1}
\end{equation}
with $F(t)$ being a Gaussian colored external force such that
\begin{equation}
\dot{F}\left( t\right) +\tau ^{-1}F\left( t\right) =\tau ^{-1}W\left(
t\right) ,  \label{S3}
\end{equation}
and
\begin{equation*}
U\left( \xi \right) =a_{1}\xi +a_{2}\xi ^{2}+a_{3}\xi ^{3}.
\end{equation*}
Here $c$ is a damping constant, $a_{j}$, $j=1,2,3,$ are some model constant
parameters, and $\tau $ is the correlation time, while $W\left( t\right) $
is a white Gaussian noise with the correlation function
\begin{equation*}
\left\langle W(t)W(t^{\prime })\right\rangle =D\delta (t-t^{\prime }).
\end{equation*}

The parameters of the model, $c,$ $a_{1},$ $a_{2},$ $a_{3},$ $D$ and $\tau $%
, are directly related to the stability of the signal: $a_{1}$ coincides
with $w^{2}$, $w=2\pi f$ being the fundamental frequency, and the parameters
$c$ and $a_{1}$ determine $DR$, see Appendix.

\begin{figure}
\begin{center}
\includegraphics[width=0.9\textwidth]{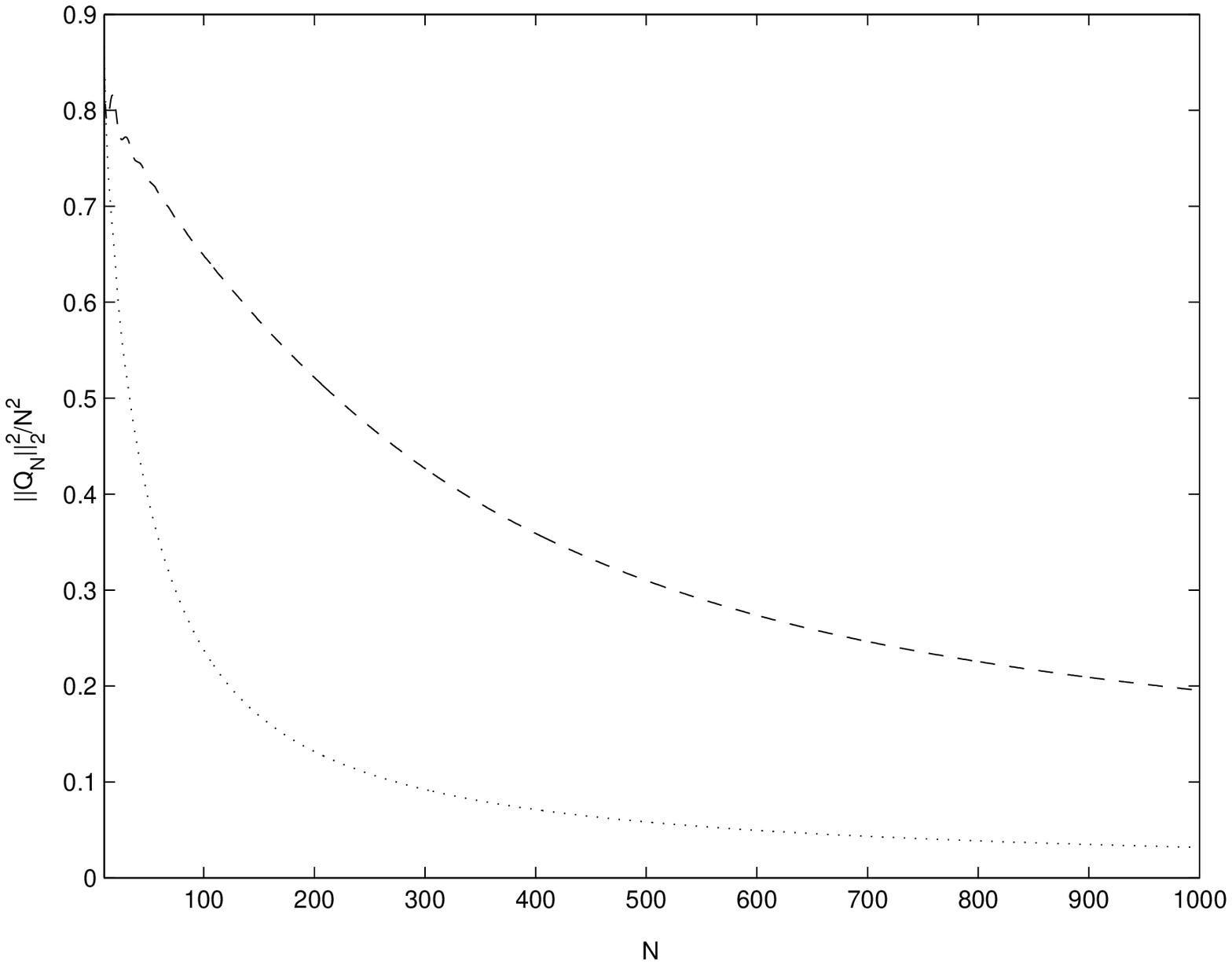}
\caption{\label{fig:5}Sequences of $\frac{1}{N^{2}}%
\left\| Q_{N}\right\| _{2}^{2}$ for the model signals, linear scale of $N$.}
\end{center}
\end{figure}

\begin{figure}
\begin{center}
\includegraphics[width=0.9\textwidth]{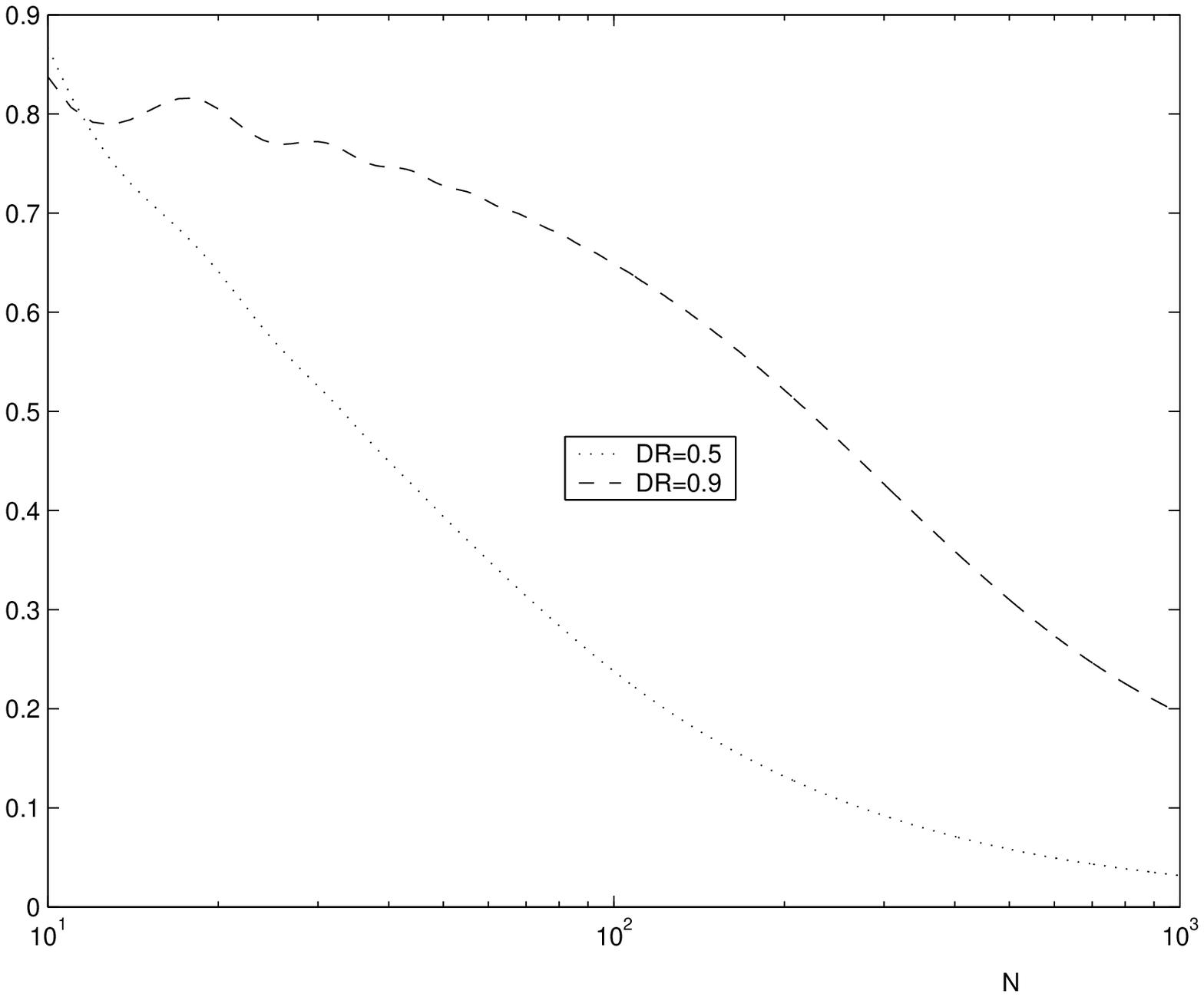}
\caption{\label{fig:6} Sequences of $\frac{1}{N^{2}}%
\left\| Q_{N}\right\| _{2}^{2}$ for the model signals, logarithmical scale
of $N$.}
\end{center}
\end{figure}

In order to compare the results obtained for the simulated signals with
those obtained previously for the real ones, we constructed the simulated
signals with the values of the parameters chosen to produce the values of $f$
and $DR$ similar to those of the signals \textit{c4\_lprm3} and \textit{%
c4\_lprm22}, and with the same sampling interval, $0.08$ $s$, and the
same number of points, 4209. In Fig.\ref{fig:5} and Fig.\ref{fig:6} the sequences of the values
of the T\"{o}plitz matrices for the models with $c=0.689$ ($DR=0.5$)
and $c=0.105 $ ($DR=0.9$) are gathered (the values of other parameters
are: $a_{1}=9$ $.87,$which corresponds to $f=0$ $.5$ $Hz$, $a_{2}=a_{3}=0,$ $%
D=500, $ $\tau =0 $ $.6$).

As expected, the behavior of these signals is close to that obtained for the
real signals, the sequences tend slower to zero for higher $DR$.

For the same simulated signals, the discontinuity points are shown in Fig.7
and Fig.8. As before, the number of segments in which the interval [-$\pi $
, $\pi $] was split was equal to 3000.

\begin{figure}
\begin{center}
\includegraphics[width=0.9\textwidth]{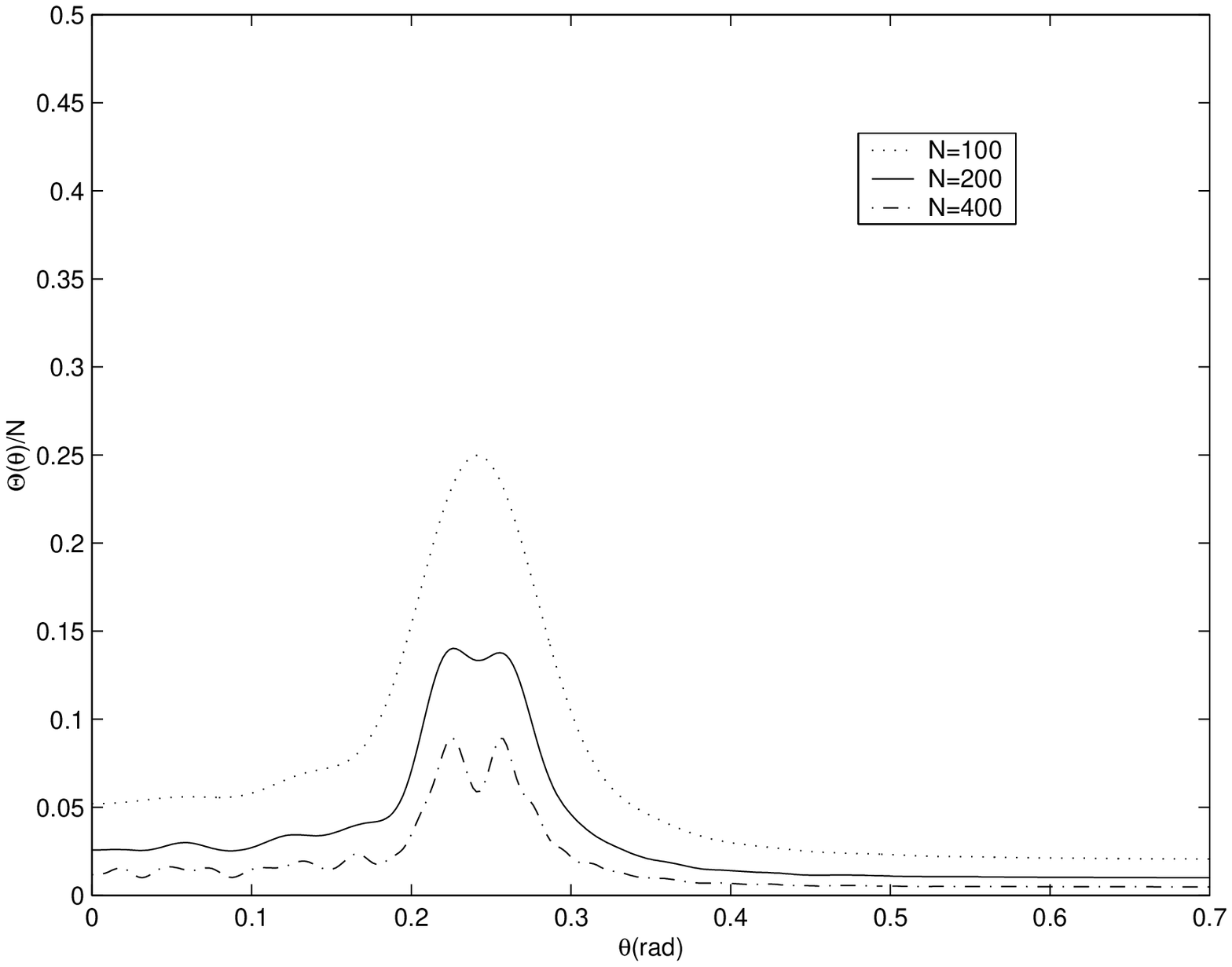}
\caption{\label{fig:7} Function $\frac{1}{N}\Theta
_{N}\left( \theta \right) $, for the model with $DR=0.5,$ $N=100,$ $N=200,$ $
N=400$.}
\end{center}
\end{figure}

\begin{figure}
\begin{center}
\includegraphics[width=0.9\textwidth]{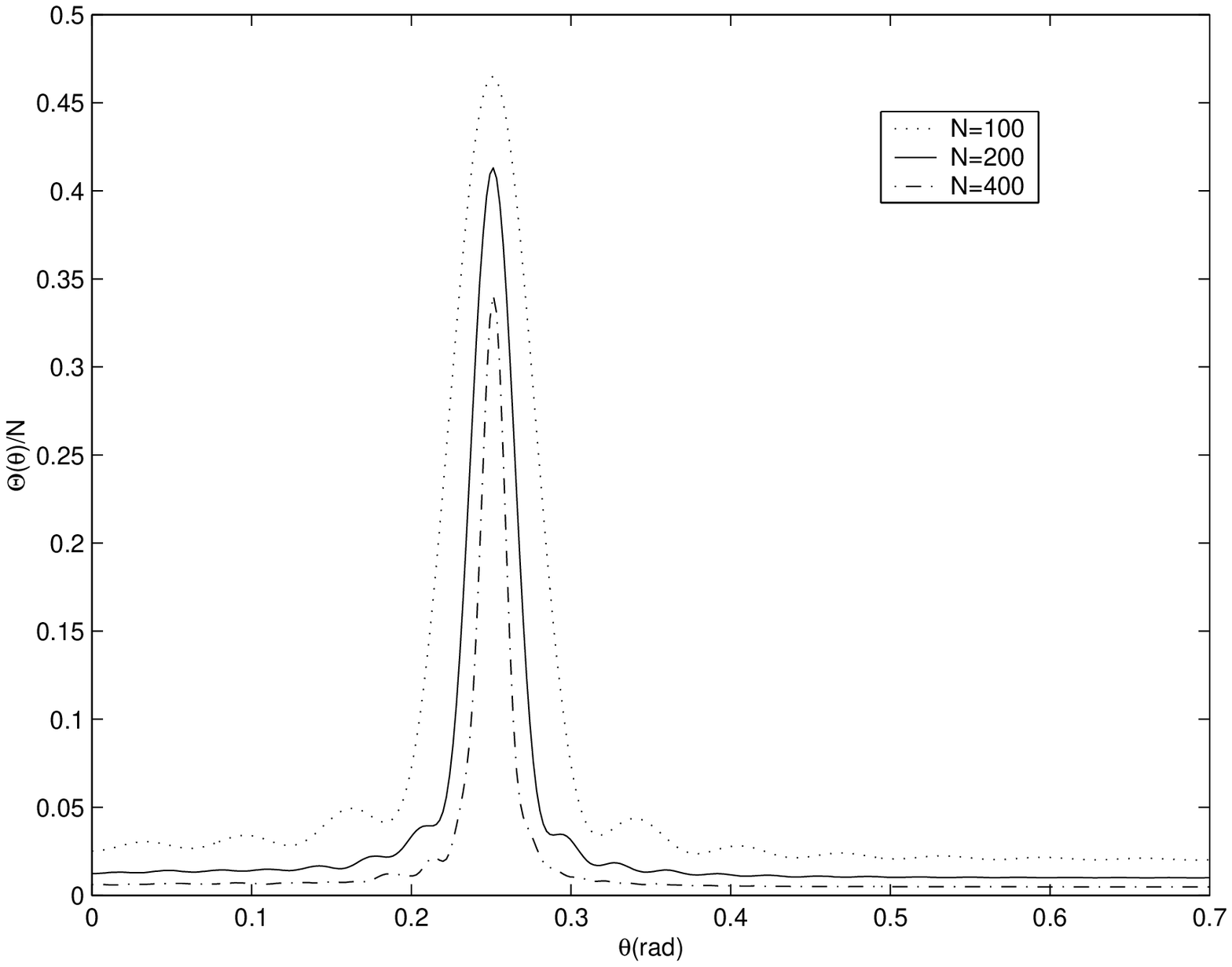}
\caption{\label{fig:8} Function $\frac{1}{N}\Theta
_{N}\left( \theta \right) $, for the model with $DR=0.9,$ $N=100,$ $N=200,$ $%
N=400.$}
\end{center}
\end{figure}

Again, for lower $DR$ the function $\frac{1}{N}\Theta _{N}\left( \theta
\right) $ tends to zero with $N$ much more rapidly. The main peak is now
obtained at $\theta =0$ $.24$ ($f=0$ $.48$ $Hz)$ for the model with $DR=0$ $%
.5$, and at $\theta =0$ $.25$ rad $(f=0$ $.50$ $Hz)$ for the model with $%
DR=0 $ $.9$.

\subsection{Appendix}

The notion of decay ratio ($DR$), a basic parameter in the analysis of
reactor stability, is deduced from the oscillatory model (\ref{S1}), but
with \cite{STVMBG}
\begin{equation}
U\left( \xi \right) =w^{2}\xi \text{ .}  \label{S2}
\end{equation}
The $DR$ is a measure of the system damping defined as the ratio between two
consecutive maxima of the signal, for the model (\ref{S1}), (\ref{S2}) it is
a constant parameter,
\begin{equation*}
DR=\exp \left\{ -\frac{2\pi c}{\sqrt{4\omega ^{2}-c^{2}}}\right\} \text{ }.
\end{equation*}

Neutronic signals are very noisy and, in general, their behavior cannot be
fitted to that of a continuous second-order system, hence the $DR$ in
reality is not a constant, and its value depends on the model used to
evaluate it \cite{VGMNPLCRS}. Nevertheless, it gives a hint on the system
stability.

\end{document}